\documentclass[a4paper,12pt]{article}
\usepackage{comment}
\usepackage{cite}
\usepackage{amsmath}
\usepackage{amssymb}
\usepackage{amsfonts}
\usepackage[T1]{fontenc}
\usepackage[utf8]{inputenc}
\usepackage{graphicx}
\usepackage{fancyhdr}
\usepackage{float}
\usepackage{xcolor}
\usepackage{authblk}
\usepackage{mathrsfs}
\usepackage{empheq}
\usepackage[hyphens]{url}
\usepackage{hyperref} 
\usepackage[]{breakurl}

\usepackage{geometry}
\begin{document}

\title{Representations of the Glaisher-Kinkelin constant deduced from integrals due to F\'eaux and Kummer}

\author[$\dagger$]{Jean-Christophe {\sc Pain}$^{1,2,}$\footnote{jean-christophe.pain@cea.fr}\\
\small
$^1$CEA, DAM, DIF, F-91297 Arpajon, France\\
$^2$Universit\'e Paris-Saclay, CEA, Laboratoire Mati\`ere en Conditions Extr\^emes,\\ 
F-91680 Bruy\`eres-le-Ch\^atel, France
}

\maketitle

\begin{abstract}
We present two integral representations of the logarithm of the Glaisher-Kinkelin constant. The calculations are based on definite integral expressions of $\log\left[\Gamma(x)\right]$, $\Gamma$ being the usual Gamma function, due respectively to F\'eaux and Kummer. The connection between the second formula and an infinite series expansion of the Glaisher-Kinkelin constant is also outlined.
\end{abstract}

\section{Introduction}\label{sec1}

The Glaisher-Kinkelin constant $A$ \cite{Finch2003} enters many sums and integrals and plays therefore an important role in number theory. For instance, it is related to the efficiency of the Euclidean algorithm and to solutions of Painlev\'e differential equations. It can be defined as
\begin{equation}
    A=\lim _{n\rightarrow \infty }{\frac {\left(2\pi \right)^{\frac {n}{2}}n^{{\frac {n^{2}}{2}}-{\frac {1}{12}}}~e^{-{\frac {3n^{2}}{4}}+{\frac {1}{12}}}}{\prod_{k=1}^{n-1}k!}}.
\end{equation}
Only a few integral representations are available in the literature (see for instance \cite{Glaisher1878,Almkvist1998,Choi1997}). We recently obtained \cite{Pain2024a}:
\begin{equation*}
\int_{0}^{\infty}\frac {(1-e^{-x/2})(x\mathrm{cotanh}\left(\frac{x}{2}\right)-2)}{x^3}\,\mathrm{d}x=3\log A-\frac {1}{3}\log 2-\frac{1}{8}
\end{equation*}
as well as
\begin{equation*}
\int_{0}^{\infty}\frac {(8-3x)e^{x}-8e^{x/2}-x}{4x^{2}e^{x}(e^{x}-1)}\,\mathrm{d}x=3\log A-\frac{7}{12}\log 2+\frac{\log\pi}{2}-1.
\end{equation*}
In sections \ref{sec2} and \ref{sec3}, we derive, using the F\'eaux and Kummer integrals (for the logarithm of the Gamma function) respectively, two integral representations for the Glaisher-Kinkelin constant. 

\section{Integral representation resulting from the F\'eaux formula}\label{sec2}

The F\'eaux formula reads \cite{Whittaker1990,Campbell1966,Gilbert1875}:
\begin{equation}\label{feaux}
    \log\left[\Gamma(x+1)\right]=\int_0^{\infty}\left[x\,e^{-t}+\frac{(1+t)^{-x-1}-(1+t)^{-1}}{\log(1+t)}\right]\frac{\mathrm{d}t}{t}.
\end{equation}
Note that there is a typographical error in Ref. \cite{Campbell1966} (Ch. V, p. 189): $(1+t)$ should be replaced by $\log(1+t)$ in the denominator. The formula is correct in Ref. \cite{Whittaker1990} (formula 12.3.4 p. 259; subsection 12.3).

Let us consider the identity \cite{Glaisher1878}:
\begin{equation}\label{gla}
    \int_0^{1/2}\log\left[\Gamma(x+1)\right]\,\mathrm{d}x=-\frac{1}{2}-\frac{7}{24}\log 2+\frac{\log\pi}{4}+\frac{3}{2}\log A, 	
\end{equation}
where $\Gamma$ represents the usual Gamma function. Using the Fubini theorem, one can write
\begin{equation*}
    \int_0^{1/2}\log\left[\Gamma(x+1)\right]\,\mathrm{d}x=\int_0^{\infty}\frac{e^{-t}}{t}~\left(\int_0^{1/2}\left[x\,e^{-t}+\frac{(1+t)^{-x-1}-(1+t)^{-1}}{\log(1+t)}\right]\,\mathrm{d}x\right)\,\mathrm{d}t,
\end{equation*}
which is equal to
\begin{equation}
    \int_0^{1/2}\log\left[\Gamma(x+1)\right]\,\mathrm{d}x=\int_0^{\infty}\left[\frac{e^{-t}}{8}-\frac{1}{(1+t)^{3/2}\log^2(1+t)}-\frac{1}{2}\frac{(\log(1+t)-2)}{(1+t)\log^2(1+t)}\right]~\frac{\mathrm{d}t}{t},
\end{equation}
or equivalently
\begin{equation}
    \int_0^{1/2}\log\left[\Gamma(x+1)\right]\,\mathrm{d}x=\int_0^{\infty}\left[\frac{e^{-t}}{8}-\frac{1}{(1+t)\log^2(1+t)}\left(\frac{1}{\sqrt{1+t}}+\frac{1}{2}(\log(1+t)-2)\right)\right]~\frac{\mathrm{d}t}{t}.
\end{equation}
Inserting the latter expression in Eq. (\ref{gla}) yields
\begin{align*}
    -\frac{1}{2}-\frac{7}{24}\log 2&+\frac{\log\pi}{4}+\frac{3}{2}\log A=\int_0^{1/2}\log\left[\Gamma(x+1)\right]\,\mathrm{d}x\nonumber\\
    =&\int_0^{\infty}\left[\frac{e^{-t}}{8}-\frac{1}{(1+t)^{3/2}\log^2(1+t)}-\frac{1}{2}\frac{(\log(1+t)-2)}{(1+t)\log^2(1+t)}\right]~\frac{\mathrm{d}t}{t}, 	
\end{align*}
leading to
\begin{align}\label{res1}
    \log A=\frac{1}{3}&+\frac{7}{36}\log 2-\frac{\log\pi}{6}\nonumber\\
    &+\frac{2}{3}\int_0^{\infty}\left[\frac{e^{-t}}{8}-\frac{1}{(1+t)^{3/2}\log^2(1+t)}-\frac{1}{2}\frac{(\log(1+t)-2)}{(1+t)\log^2(1+t)}\right]~\frac{\mathrm{d}t}{t},  	
\end{align}
which is the first main result of the present work. 

\section{Integral representation resulting from the Kummer formula}\label{sec3}

The Kummer formula reads \cite{Erdelyi1981,Kummer1847}:
\begin{equation*}\label{kum}
    \log\left[\Gamma(x)\right]=\frac{\log\pi}{2}-\frac{1}{2}\log\left[\sin(\pi x)\right]+\frac{1}{2}\int_0^{\infty}\left[\frac{\sinh\left[\left(\frac{1}{2}-x\right)t\right]}{\sinh\left(\frac{t}{2}\right)}-(1-2x)\,e^{-t}\right]\,\frac{\mathrm{d}t}{t}.
\end{equation*}
Let us start from the identity:
\begin{equation}\label{gla2}
\log A=\frac{2}{3}\int_0^{1/2}\log\left[\Gamma(x)\right]\,\mathrm{d}x-\frac{5}{36}\log 2-\frac{\log\pi}{6},
\end{equation}
which is a particular case of the expression \cite{Alexejewsky1894,Barnes1899,Adamchik2003}:
\begin{equation*}
\int_{0}^{z}\log \Gamma (x)\,dx={\frac {z(1-z)}{2}}+{\frac {z}{2}}\log 2\pi +z\log \Gamma (z)-\log G(1+z),
\end{equation*}
where $G$ is the Barnes $G$-function. Using the Fubini theorem, one can write
\begin{equation}\label{inter}
    \int_0^{1/2}\log\left[\Gamma(x)\right]\,\mathrm{d}x=\frac{\log(2\pi)}{4}+\int_0^{\infty}\left(\int_0^{1/2}\left[\frac{\sinh\left[\left(\frac{1}{2}-x\right)t\right]}{\sinh\left(\frac{t}{2}\right)}-(1-2x)\,e^{-t}\right]\,\mathrm{d}x\right)\,\frac{\mathrm{d}t}{t},
\end{equation}
where we have used
\begin{equation}
    \int_0^{1/2}\log\left[\sin(\pi x)\right]\,\mathrm{d}x=-\frac{\log 2}{2}.
\end{equation}
Equation (\ref{inter}) becomes
\begin{equation*}
    \int_0^{1/2}\log\left[\Gamma(x)\right]\,\mathrm{d}x=\frac{\log(2\pi)}{4}+\frac{1}{2}\int_0^{\infty}\left[\frac{\tanh\left(\frac{t}{4}\right)}{t}-\frac{e^{-t}}{4}-\right]\,\frac{\mathrm{d}t}{t}.
\end{equation*}
Inserting the latter expression in Eq. (\ref{gla2}) yields
\begin{equation*}
    \log A=-\frac{5}{36}\log 2-\frac{\log\pi}{6}+\frac{\log (2\pi)}{6}+\frac{1}{3}\int_0^{\infty}\left[\frac{\tanh\left(\frac{t}{4}\right)}{t}-\frac{e^{-t}}{4}-\right]\,\frac{\mathrm{d}t}{t}, 	
\end{equation*}
leading to
\begin{equation}\label{res2}
    \log A=\frac{\log 2}{36}+\frac{1}{3}\int_0^{\infty}\left[\frac{\tanh\left(\frac{t}{4}\right)}{t}-\frac{e^{-t}}{4}-\right]\,\mathrm{d}t, 	
\end{equation}
which is the second main result of the present work. The Kummer integral can also help obtaining an infinite series expansion of the Glaisher-Kinkelin constant (see Appendix).

\section{Conclusion}

In this article, we presented two integral representations of the logarithm of the Glaisher-Kinkelin constant. Both are based on a definite integral representation involving the logarithm of the Gamma function, which can itself be expressed by an improper integral. In that perspective, our formulas are obtained using two results, one due to F\'eaux, and the other due to Kummer. It is hoped that the new integral representations will provide new insights on the properties of the Glaisher-Kinkelin constant.

\section*{Appendix: Infinite series for the Glaisher-Kinkelin constant}

By expanding $\sin\left[(x-1/2)t\right]$ and $(1-2x)$ of the integral representation (\ref{kum}) in Fourier sine series, one gets
\begin{equation}\label{fou}
    \log\left[\Gamma(x)\right]=\frac{\log\pi}{2}-\frac{1}{2}\log\left[\sin(\pi x)\right]+2\sum_{n=1}^{\infty}a_n\,\sin(2\pi nx),
\end{equation}
with
\begin{equation}
    a_n=\int_0^{\infty}\left\{\frac{2n\pi}{t^2+4n^2\pi^2}-\frac{e^{-t}}{2n\pi}\right\}\frac{\mathrm{d}t}{t}.
\end{equation}
Furthermore, using the Dirichlet formula,
\begin{equation}
    \gamma=\int_0^{\infty}\left(\frac{1}{1+t}-e^{-t}\right)\frac{\mathrm{d}t}{t},
\end{equation}
one obtains \cite{Whittaker1990,Farhi2013,Pain2024c}:
\begin{equation}\label{an}
a_n=\frac{1}{2n\pi}(\gamma+\log(2\pi)+\log n).
\end{equation}
Inserting Eq. (\ref{fou}) into Eq. (\ref{gla2}) provides in infinite series representation of the logarithm of the Glaisher-Kinkelin function:
\begin{equation*}
    \log A=\frac{\log 2}{36}+\frac{2}{3\pi}\sum_{n=1}^{\infty}\left[1-(-1)^n\right]\frac{a_n}{n},
\end{equation*}
i.e.,
\begin{equation*}
    \log A=\frac{\log 2}{36}+\frac{4}{3\pi}\sum_{n=1}^{\infty}\frac{a_{2n+1}}{2n+1},
\end{equation*}
and replacing $a_{2n+1}$ by its expression (\ref{an}), one gets
\begin{equation*}
    \log A=\frac{\log 2}{36}+\frac{2(\gamma+\log(2\pi))}{3\pi^2}\sum_{n=1}^{\infty}\frac{1}{(2n+1)^2}+\frac{2}{3\pi^2}\sum_{n=0}^{\infty}\frac{\log(2n+1)}{(2n+1)^2},
\end{equation*}
leading, using
\begin{equation}
    \sum_{n=0}^{\infty}\frac{1}{(2n+1)^2}=\frac{\pi^2}{8},
\end{equation}
to the final result
\begin{equation*}
    \log A=\frac{\log 2}{36}+\frac{\gamma+\log(2\pi)}{12}+\frac{2}{3\pi^2}\sum_{n=0}^{\infty}\frac{\log(2n+1)}{(2n+1)^2}.
\end{equation*}
Note that Hasse obtained, using a series or the Riemann zeta function \cite{Guillera2008}:
\begin{equation*}
\log A=\frac{1}{8}-\frac{1}{2}\sum_{n=0}^{\infty}\frac{1}{n+1}\sum_{k=0}^{n}(-1)^{k}\binom{n}{k}(k+1)^{2}\log(k+1).
\end{equation*}

\end{document}